\documentclass[a4paper,10pt]{amsart}

\usepackage{acatta}
\usepackage{enumitem}

\newcommand{\blankSpace}{\cdot}

\title{A note on the automorphism group of Kodaira surfaces}

\author{Andrea Cattaneo}
\address{Universit\`a di Parma\\
Dipartimento di Scienze Matematiche, Fisiche e Informatiche\\
Unit\`a di Matematica e Informatica\\
Parco Area delle Scienze 53/A, 43124, Parma, Italy}
\email{andrea.cattaneo@unipr.it}


\begin{document}

\begin{abstract}
We give an explicit description of the automorphism groups of primary and secondary Kodaira surfaces. We first prove that every automorphism of a Kodaira surface lifts to an affine transformation of the universal cover $\IC^2$, and characterize precisely the affine transformations arising in this way for primary Kodaira surfaces. This allows us to compute explicitly the normalizer of the fundamental group and, consequently, to determine both the connected component of the identity and the group of connected components of the automorphism group. Passing to the canonical cover, we then obtain an analogous explicit description for secondary Kodaira surfaces, including a complete determination of the possible groups of connected components.
\end{abstract}

\maketitle

\tableofcontents

\section{Introduction}

The study of automorphism groups of compact complex manifolds is a classical topic in complex geometry and has attracted renewed interest in recent years. Besides their intrinsic geometric interest, automorphism groups often reflect subtle properties of the underlying manifolds and play a central role in problems concerning moduli, dynamics and birational geometry.

Among compact complex surfaces of Kodaira dimension zero, Kodaira surfaces occupy a distinguished position, being the only non-K\"ahler examples. They also admit one of the most concrete presentations among compact complex surfaces, being realized as quotients of the complex affine plane by properly discontinuous groups of affine transformations. This naturally suggests studying their automorphisms through the action on the universal covering space: every automorphism lifts to an automorphism of the plane normalizing the corresponding deck transformation group. The guiding idea of this paper is that this classical presentation already contains enough information to recover the full structure of the automorphism groups of both primary and secondary Kodaira surfaces. In this paper we exploit this presentation in order to obtain an explicit computation of these groups.

We begin by considering primary Kodaira surfaces. Our first main result shows that every lift of an automorphism to the universal cover is necessarily affine (Lemma \ref{lemma: lift 1}). This allows us to compute explicitly the normalizer of the fundamental group inside $\Aut(\IC^2)$ (Proposition \ref{prop: lift 2}). From this computation we derive an explicit description of the full automorphism group: we identify its connected component with the elliptic fibre of the Albanese fibration and determine the group of connected components by describing it as an extension of the finite cyclic group $\mu_n$ by $(\IZ \oplus \IZ/m\IZ)^2$, leading to the structural description summarized in Theorem \ref{thm: aut primary kodaira}.

We then apply the same philosophy to secondary Kodaira surfaces. Since every automorphism lifts to the canonical cover, which is a primary Kodaira surface, the previous computation can be combined with a description of the normalizer of the covering group. This yields an analogous explicit description of the automorphism group of every secondary Kodaira surface (Theorem \ref{thm: list of groups}), including a complete determination of the possible groups of connected components.

Finally, we compare our description with that obtained by Fujimoto and Na\-ka\-ya\-ma for surjective endomorphisms of Kodaira surfaces. We point out a slight inaccuracy in the proof of their description of the lifted maps and state the corresponding corrected result, thereby clarifying the consequences for subsequent work relying on that proposition.

The paper is organized as follows. In Section \ref{sect: Kod surf} we briefly recall Kodaira's original description of primary Kodaira surfaces and establish the explicit computation of the normalizer, leading to the description of their automorphism groups. Section \ref{sect: secondary Kodaira} extends the same approach to secondary Kodaira surfaces through the canonical covering. Finally, Section \ref{sect: comparison with FN} discusses the comparison with Fujimoto–Nakayama and provides the corrected statement.

\begin{ack}
The author is a member of the GNSAGA of INdAM and was supported by the project PRIN2022 ``Real and Complex Manifolds: Geometry and Holomorphic Dynamics'' (Project code: 2022AP8HZ9). This research was granted by University of Parma through the action ``Bando di Ateneo 2025 per la ricerca''. He wants to warmly acknowledge professor Fabrizio Catanese for his useful suggestions and encouragement to put this paper in the present form, and for having pointed out his survey \cite{Catanese} and Proposition \ref{prop: exact sequence}.
\end{ack}

\section{Primary Kodaira surfaces}\label{sect: Kod surf}

In \cite[$\S$6]{Kodaira1} Kodaira studied the compact complex surfaces with trivial can\-on\-i\-cal bundle, proving the following theorem.

\begin{thm}[{\cite[Theorem 19]{Kodaira1}}]
Let $S$ be a surface with trivial canonical bundle. Then $S$ is a K3 surface, a complex torus or an elliptic surface of the form $\IC^2 / G$ where $G$ is a properly discontinuous group of affine transformations without fixed points of $\IC^2$ which leave invariant the $2$-form $dz \wedge d\zeta$ (being $(z, \zeta)$ the standard coordinates). The first homology group of the elliptic surface $\IC^2 / G$ is $H_1(\IC^2 / G, \IZ) \simeq \IZ^3 \oplus \IZ / m\IZ$.
\end{thm}

From this theorem, a \emph{primary Kodaira surface} is a compact complex surface $X$ with trivial canonical bundle and odd first Betti number. It is known (see, e.g., \cite[Teorema 52]{Catanese}) that the universal cover of a primary Kodaira surface is $\IC^2$ and that its fundamental group $G$ is a central extension of $\IZ^3 \oplus \IZ / m\IZ$ by $\IZ$, which can be generated by the four transformations
\begin{equation}\label{eq: generators of G}
\begin{array}{ll}
g_1(z, \zeta) = (z, \zeta + \beta_1), & g_2(z, \zeta) = (z, \zeta + \beta_2),\\
g_3(z, \zeta) = (z + \alpha_3, \zeta + \bar{\alpha}_3 z + \beta_3), & g_4(z, \zeta) = (z + \alpha_4, \zeta + \bar{\alpha}_4 z + \beta_4),
\end{array}
\end{equation}
where $\beta_1, \beta_2$ and $\alpha_3, \alpha_4$ are two pairs of complex numbers which are linearly independent over $\IR$ and satisfy the relation
\begin{equation}\label{eq: relation defying m}
\bar{\alpha}_3 \alpha_4 - \alpha_3 \bar{\alpha}_4 = m \beta_2
\end{equation}
for some $m \in \IZ_{> 0}$. This relation is equivalent to the relation
\[[g_3, g_4] = g_3 \circ g_4 \circ g_3^{-1} \circ g_4^{-1} = g_2^m.\]

\begin{rem}\label{rem: Albanese}
The projection on the first factor $\IC^2 \longrightarrow \IC$ is equivariant with respect to the natural action of $G \leq \Aut(\IC^2)$ on $\IC^2$ and the action given by translations via elements in $\Lambda_B = \IZ \cdot \alpha_3 \oplus \IZ \cdot \alpha_4$ on $\IC$. As a consequence, it induces on the quotients the map
\[a: X \longrightarrow B = \IC / \Lambda_B,\]
where $B$ is a compact complex curve of genus $1$. It turns out that this map is the Albanese map of $X$ and that its fibres are isomorphic to the compact complex torus $E = \IC / \Lambda_E$, where $\Lambda_E = \IZ \cdot \beta_1 \oplus \IZ \cdot \beta_2$.
\end{rem}

\begin{rem}\label{rem: elements in G}
It is possible to describe $G$ in another way, as a central extension of $\Lambda_B$ by $\Lambda_E$. Accordingly, every element $g \in G$ can be written in a unique way as
\[g(z, \zeta) = (z + \alpha, \zeta + \bar{\alpha} z + \beta + P(\alpha))\]
for suitable $\alpha \in \Lambda_B$ and $\beta \in \Lambda_E$, and where for every $\alpha = m_3 \alpha_3 + m_4 \alpha_4 \in \Lambda_B$ we define
\[P(\alpha) = \frac{1}{2} m_3 (m_3 - 1) |\alpha_3|^2 + m_3 m_4 \bar{\alpha}_3 \alpha_4 + \frac{1}{2} m_4 (m_4 - 1) |\alpha_4|^2 + m_3 \beta_3 + m_4 \beta_4.\]
Observe that then
\begin{equation}\label{eq: relation P}
P(\alpha) + P(-\alpha) = |\alpha|^2 + m_3 m_4 m \beta_2.
\end{equation}
\end{rem}

Since $\IC^2$ is the universal covering of $X$, the group $\Aut(X)$ is given by
\[\Aut(X) = N_{\Aut(\IC^2)}(G) / G\]
where $N_{\Aut(\IC^2)}(G) = \set{f \in \Aut(\IC^2) \st f \circ G = G \circ f}$ is the normalizer of $G$ in $\Aut(\IC^2)$.

\subsection{Computation of the normalizer}\label{sect: normalizer}

Through this section, $X$ will denote a primary Kodaira surface. Let $f \in \Aut(X)$ and denote by $F: \IC^2 \longrightarrow \IC^2$ a lift of $f$ to the universal cover. Then $F(z, \zeta) = (F_1(z, \zeta), F_2(z, \zeta))$ is a holomorphic map with the property that for every $g \in G$ there exists $\gamma = \gamma(g) \in G$ such that
\begin{equation}\label{eq: normalizing}
F \circ g = \gamma \circ F.
\end{equation}

\begin{lemma}\label{lemma: lift 1}
Let $X$ be a primary Kodaira surface and let $f \in \Aut(X)$. Let $F: \IC^2 \longrightarrow \IC^2$ be a lift of $f$ to the universal cover, then $F$ is an affine transformation of $\IC^2$ of the form
\[F(z, \zeta) = (vz + \delta, |v|^2 \zeta + \mu z + \sigma)\]
for suitable $v, \delta, \mu, \sigma \in \IC$ with $v \neq 0$.
\end{lemma}
\begin{proof}
Since the generic map $g \in G$ has the form (compare with Remark \ref{rem: elements in G})
\[g(z, \zeta) = (z + \alpha_g, \zeta + \bar{\alpha}_g z + \beta_g + P(\alpha_g)), \qquad (\beta_g, \alpha_g) \in \Lambda_E \times \Lambda_B,\]
taking the Jacobian matrix of \eqref{eq: normalizing} gives
\begin{equation}\label{eq: jacobian matrix}
\begin{split}
&\left(
\begin{array}{cc}
\frac{\partial F_1}{\partial z}(g(z, \zeta)) &
\frac{\partial F_1}{\partial \zeta}(g(z, \zeta)) \\
\frac{\partial F_2}{\partial z}(g(z, \zeta)) &
\frac{\partial F_2}{\partial \zeta}(g(z, \zeta))
\end{array}
\right)
\left(
\begin{array}{cc}
1 & 0 \\
\bar{\alpha}_g & 1
\end{array}
\right)=
\\
&=
\left(
\begin{array}{cc}
1 & 0 \\
\bar{\alpha}_{\gamma(g)} & 1
\end{array}
\right)
\left(
\begin{array}{cc}
\frac{\partial F_1}{\partial z}(z, \zeta) &
\frac{\partial F_1}{\partial \zeta}(z, \zeta) \\
\frac{\partial F_2}{\partial z}(z, \zeta) &
\frac{\partial F_2}{\partial \zeta}(z, \zeta)
\end{array}
\right).
\end{split}
\end{equation}
We proceed by examining the four entries of this matrix.

First, we focus on the entry $(1, 2)$ of \eqref{eq: jacobian matrix}. Explicitly, we find that
\[\frac{\partial F_1}{\partial \zeta}(g(z, \zeta)) = \frac{\partial F_1}{\partial \zeta}(z, \zeta), \qquad \text{for every } g \in G,\]
which means that $\frac{\partial F_1}{\partial \zeta}$ is a well defined function on $X$.
Since $X$ is compact this function must be constant, hence
\[\frac{\partial F_1}{\partial \zeta}(z, \zeta) = k \in \IC\]
and so
\begin{equation}\label{eq: 1,2}
F_1(z, \zeta) = k\zeta + H_1(z)
\end{equation}
for a suitable holomorphic function $H_1$.

Then we address to the entry $(1, 1)$ of \eqref{eq: jacobian matrix}. Using \eqref{eq: 1,2} we see that $H_1'(z + \alpha_g) + \bar{\alpha}_g k = H_1'(z)$ and so
\[H_1''(z + \alpha_g) = H_1''(z).\]
Hence $H_1''$ descends to a well defined function on the compact curve $B$ and so it is constant, say $H_1''(z) = u \in \IC$. But then $H_1(z) = \frac{1}{2}uz^2 + vz + \delta$ for some $v, \delta \in \IC$ and so
\[F_1(z, \zeta) = \frac{1}{2}uz^2 + vz + k\zeta + \delta.\]
As a consequence of \eqref{eq: normalizing} we deduce that for every $g \in G$ the following conditions hold:
\[\alpha_g u + \bar{\alpha}_g k = 0, \qquad \alpha_{\gamma(g)} = \frac{1}{2} u \alpha_g^2 + v \alpha_g + k \beta_g.\]
The first relation holds in particular for $g_3$ and $g_4$ (where $\alpha_{g_i} = \alpha_i$ for $i = 3, 4$ according to our notation): from this we deduce that $u = k = 0$, and so
\begin{equation}\label{eq: 1,1}
F_1(z, \zeta) = v z + \delta \qquad \text{and} \qquad \alpha_{\gamma(g)} = v \alpha_g.
\end{equation}

Next, look at the the entry $(2, 2)$ of \eqref{eq: jacobian matrix}. As $F_1$ does not depend on $\zeta$ we see that
\[\frac{\partial F_2}{\partial \zeta}(g(z, \zeta)) = \frac{\partial F_2}{\partial \zeta}(z, \zeta), \qquad \text{for every } g \in G.\]
So $\frac{\partial F_2}{\partial \zeta}$ descends to a function on $X$, which must be constant: $\frac{\partial F_2}{\partial \zeta}(z, \zeta) = \rho \in \IC$.

Finally, concentrate on the entry $(2, 1)$ of \eqref{eq: jacobian matrix}. Explicitly we see that
\[\frac{\partial F_2}{\partial z}(g(z, \zeta)) + \bar{\alpha}_g \rho = \bar{\alpha}_g v + \frac{\partial F_2}{\partial z}(z, \zeta),\]
from which we deduce that
\[\frac{\partial^2 F_2}{\partial z^2}(g(z, \zeta)) = \frac{\partial^2 F_2}{\partial z^2}(z, \zeta).\]
This means that $\frac{\partial^2 F_2}{\partial z^2}(z, \zeta)$ descends to a function on $X$, which is then constant. Hence
\[F_2(z, \zeta) = \frac{1}{2} \lambda z^2 + \mu z + \rho \zeta + \sigma\]
for suitable $\lambda, \mu, \sigma \in \IC$. Using \eqref{eq: normalizing} and \eqref{eq: 1,1} we deduce that for every $g \in G$ the following conditions hold:
\[\alpha_g \lambda + \bar{\alpha}_g \rho = |v|^2 \bar{\alpha}_g, \qquad \beta_{\gamma(g)} + P(\alpha_{\gamma(g)}) = \frac{1}{2} \lambda \alpha_g^2 + \mu \alpha_g + \rho \beta_g - \bar{v} \bar{\alpha}_g \delta\]
The first one holds in particular for $g_3$ and $g_4$: from this we deduce that $\lambda = 0$ and $\rho = |v|^2$, and so
\begin{equation}\label{eq: 2,1}
F_2(z, \zeta) = \mu z + |v|^2 \zeta + \sigma \qquad \text{and} \qquad \beta_{\gamma(g)} + P(\alpha_{\gamma(g)}) = \mu \alpha_g + |v|^2 \beta_g - \bar{v} \bar{\alpha}_g \delta.
\end{equation}

Observe that $v \in \IC$ can not be zero, otherwise $f$ could not be onto, which means that $F$ is invertible and we have
\begin{equation}\label{eq: F inverse}
F^{-1}(z, \zeta) = \left( \frac{1}{v}z - \frac{\delta}{v}, \frac{1}{|v|^2} \zeta - \frac{\mu}{v |v|^2} z + \frac{\delta \mu - v \sigma}{v |v|^2} \right).
\end{equation}
\end{proof}

\begin{rem}\label{rem: FG subset GF}
Let $f$ and $F$ be as in Lemma \ref{lemma: lift 1}, and consider $g(z, \zeta) = (z + \alpha_g, \zeta + \bar{\alpha}_g z + \beta_g + P(\alpha_g)) \in G$ for $\beta_g \in \Lambda_E$ and $\alpha_g \in \Lambda_B$ as in Remark \ref{rem: elements in G}. Then
\[\gamma(z, \zeta) = (F \circ g \circ F^{-1})(z, \zeta) = (z + v \alpha_g, \zeta + \bar{v} \bar{\alpha}_g z + |v|^2(\beta_g + P(\alpha_g)) - \bar{v} \bar{\alpha}_g \delta + \mu \alpha_g)\]
and so:
\begin{enumerate}
\item\label{item: Lambda_E} for $g$ of the form $g(z, \zeta) = (z, \zeta + \beta_g)$ with $\beta_g \in \Lambda_E$ we deduce that $|v|^2 \beta_g \in \Lambda_E$, i.e., $|v|^2 \cdot \Lambda_E \subseteq \Lambda_E$;
\item\label{item: Lambda_B} $\alpha_\gamma = v \alpha_g \in \Lambda_B$ for every $\alpha_g \in \Lambda_B$, i.e., $v \cdot \Lambda_B \subseteq \Lambda_B$;
\item\label{item: condition} there exists $\beta_\gamma \in \Lambda_E$ such that
\[|v|^2(\beta_g + P(\alpha_g)) - \bar{v} \bar{\alpha}_g \delta + \mu \alpha_g = \beta_\gamma + P(\alpha_\gamma) = \beta_\gamma + P(v \alpha_g).\]
\end{enumerate}
These properties are equivalent to the fact that $F \circ G \circ F^{-1} \subseteq G$.
\end{rem}

The previous discussion allows us to characterize the normalizer $N_{\Aut(\IC^2)}(G)$ explicitly.

\begin{prop}\label{prop: lift 2}
Let $X$ be a primary Kodaira surface. Then the normalizer $N_{\Aut(\IC^2)}(G)$ consists of the affine transformations
\begin{equation}\label{eq: lift final}
F(z, \zeta) = (vz + \delta, \zeta + \mu z + \sigma)
\end{equation}
for suitable $v, \delta, \mu, \sigma \in \IC$ such that
\begin{enumerate}[label = (\roman*), ref = \roman*]
\item\label{item: cond norm 1} $v \in \mu_n$, where $\Aut(B) = B \rtimes \mu_n$. In particular, $v \cdot \Lambda_B = \Lambda_B$ and $n \in \set{2, 4, 6}$;
\item\label{item: cond norm 2} $\bar{\alpha}_g \delta - \alpha_g \bar{v} \mu + P(\alpha_g) - P(\bar{v} \alpha_g) \in \Lambda_E$ for every $\alpha_g \in \Lambda_B$.
\end{enumerate}
\end{prop}
\begin{proof}
By Lemma \ref{lemma: lift 1}, a holomorphic map $F$ which lifts an automorphism $f \in \Aut(X)$ is of the form
\[F(z, \zeta) = (vz + \delta, |v|^2 \zeta + \mu z + \sigma)\]
for suitable $v, \delta, \mu, \sigma \in \IC$ with $v \neq 0$ and satisfying the three requirements of Remark \ref{rem: FG subset GF}. As $F^{-1}$ lifts $f^{-1} \in \Aut(X)$, the same must hold also for $F^{-1}$, so using the explicit expression given in \eqref{eq: F inverse} we have that
\begin{enumerate}[label = (\arabic*'), ref = \arabic*']
\item\label{item': Lambda_E} $\frac{1}{|v|^2} \cdot \Lambda_E \subseteq \Lambda_E$, i.e., $\Lambda_E \subseteq |v|^2 \cdot \Lambda_E$;
\item\label{item': Lambda_B} $\frac{1}{v} \cdot \Lambda_B \subseteq \Lambda_B$, i.e., $\Lambda_B \subseteq v \cdot \Lambda_B$;
\item\label{item': condition} $-\alpha_g \frac{\mu}{v |v|^2} + \bar{\alpha}_g \frac{\delta}{|v|^2} + \frac{1}{|v|^2} (\beta_g + P(\alpha_g)) - P\left( \frac{1}{v}\alpha_g \right) \in \Lambda_E$ for every $\beta_g \in \Lambda_E$ and $\alpha_g \in \Lambda_B$.
\end{enumerate}
Recall that these last conditions mean that $F^{-1} \circ G \circ F \subseteq G$, hence we have that $F$ normalizes $G$.

The conditions \eqref{item: Lambda_B} and \eqref{item': Lambda_B} together imply that $v \cdot \Lambda_B = \Lambda_B$, hence \eqref{item: cond norm 1}. Then the expression of $F$ is as in the statement, and \eqref{item: Lambda_E} and \eqref{item': Lambda_E} are automatically satisfied. Finally, conditions \eqref{item: condition} and \eqref{item': condition} become
\[\left\{ \begin{array}{l}
- \bar{\alpha}_g \bar{v} \delta + \alpha_g \mu + P(\alpha_g) - P(v \alpha_g) \in \Lambda_E\\
\bar{\alpha}_g \delta - \alpha_g \bar{v} \mu + P(\alpha_g) - P(\bar{v} \alpha_g) \in \Lambda_E
\end{array} \right.\]
for every $\alpha_g \in \Lambda_B$. It is enough to consider the second condition. Indeed, applying it to $-v\alpha_g \in \Lambda_B$, and using \eqref{eq: relation P}, one obtains the first condition, which is therefore redundant.

\end{proof}

The Albanese map of a primar y Kodaira surface $a: X \longrightarrow B$ (see Remark \ref{rem: Albanese}) gives $X$ the structure of a principal $E$-bundle. As a consequence we have $E \subseteq \Aut(X)$ naturally, and in fact $E \simeq TG/G \simeq T/(G \cap T) \simeq \IC^2 / \Lambda_E$ with
\[T = \set{F(z, \zeta) = (z, \zeta + \sigma) \st \sigma \in \IC^2} \leq N_{\Aut(\IC^2)}(G).\]
Moreover, $E = \Aut^0(X)$ is the connected component of the identity in $\Aut(X)$ (see \cite[Remark 2]{Borcea}).

\begin{prop}\label{prop: exact sequence}
Let $X$ be a primary Kodaira surface. There is a short exact sequence
\[1 \longrightarrow K \longrightarrow \frac{\Aut(X)}{\Aut^0(X)} \longrightarrow \mu_n \longrightarrow 1\]
where $\mu_n$ is the group of $n^{\text{th}}$-roots of unity and
\begin{equation}\label{eq: K}
K \simeq \frac{\Lambda_E^2}{(\IZ \cdot (m \beta_2))^2} \simeq \left( \IZ \oplus \frac{\IZ}{m \IZ} \right)^2.
\end{equation}
\end{prop}
\begin{proof}
Let us denote $F_{(v, \delta, \mu, \sigma)}$ the automorphism of $\IC^2$ defined by
\[F_{(v, \delta, \mu, \sigma)}(z, \zeta) = (vz + \delta, \zeta + \mu z + \sigma).\]
The natural map
\[\begin{array}{rccc}
\varphi: & N_{\Aut(\IC^2)}(G) & \longrightarrow & \mu_n\\
 & F_{(v, \delta, \mu, \sigma)} & \longmapsto & v
\end{array}\]
is a surjective homomorphism with kernel
\[\ker(\varphi) = \set{F_{(1, \delta, \mu, \sigma)} \st \sigma \in \IC, \, -\bar{\alpha}_g \delta + \alpha_g \mu \in \Lambda_E \, \text{ for all } \alpha_g \in \Lambda_B}.\]
We observe that $G, T \subseteq \ker(\varphi)$ and so $\varphi$ induces the homomorphisms
\[\tilde{\varphi}: \Aut(X) \simeq \frac{N_{\Aut(\IC^2)}(G)}{G} \longrightarrow \mu_n, \qquad \hat{\varphi}: \frac{\Aut(X)}{\Aut^0(X)} \simeq \frac{N_{\Aut(\IC^2)}(G)}{TG} \longrightarrow \mu_n.\]
It follows that
\[K = \ker(\hat{\varphi}) \simeq \frac{\ker(\varphi)}{TG} \simeq \frac{\ker(\varphi)/T}{G/(G \cap T)}.\]
Now, the map
\[\begin{array}{rcl}
\set{(\delta, \mu) \in \IC^2 \st -\bar{\alpha}_g \delta + \alpha_g \mu \in \Lambda_E \text{ for all } \alpha_g \in \Lambda_B} & \longrightarrow & \ker(\varphi)/T\\
(\delta, \mu) & \longmapsto & T \cdot F_{(1, \delta, mu, 0)}
\end{array}\]
defines a group isomorphism, which identifies $G/(G \cap T) = G/\Lambda_E \simeq \Lambda_B$ with the subgroup corresponding to $(\delta, \mu) = (\alpha, \bar{\alpha})$ for $\alpha \in \Lambda_B$. Denote
\[A = \left( \begin{array}{cc}
-\bar{\alpha}_3 & \alpha_3\\
-\bar{\alpha}_4 & \alpha_4
\end{array} \right),\]
then $\ker(\varphi)/T \simeq A^{-1} \cdot \Lambda_E^2 \simeq \Hom(\Lambda_B, \Lambda_E) \simeq \Lambda_E^2$ and $\Lambda_B$ corresponds to the subgroup
\[\IZ \cdot \left( \begin{array}{c}
0\\
m \beta_2
\end{array} \right) \oplus \IZ \cdot \left( \begin{array}{c}
-m \beta_2\\
0
\end{array} \right) \simeq (\IZ \cdot (m \beta_2))^2.\]
Hence
\[K \simeq \frac{\Lambda_E^2}{(\IZ \cdot (m \beta_2))^2} \simeq \left( \IZ \oplus \frac{\IZ}{m \IZ} \right)^2.\]
\end{proof}

\begin{rem}\label{rem: action}
It follows from Proposition \ref{prop: exact sequence} that $\mu_n$ acts on $K$ as follows: $v \in \mu_n$ acts on the element of $K$ represented by $\left( \begin{array}{c} y_1\\ y_2 \end{array} \right) \in \Lambda_E^2$ sending it to the element of $K$ represented by $V \cdot \left( \begin{array}{c} y_1 \\ y_2 \end{array} \right)$, with
\[V = A \cdot \left( \begin{array}{cc} v & 0\\ 0 & \bar{v} \end{array} \right) \cdot A^{-1}.\]
\end{rem}

We summarize the previous results in the following theorem, which describes the abstract structure of the automorphism group of a primary Kodaira surface.

\begin{thm}\label{thm: aut primary kodaira}
Let $X$ be a primary Kodaira surface. Then $\Aut^0(X) \simeq E$, where $E$ is the fibre of the Albanese map $X \longrightarrow B$. Moreover, there is a central extension
\[1 \longrightarrow E \longrightarrow \Aut(X) \longrightarrow K \rtimes \mu_n \longrightarrow 1,\]
where
\[K \cong \left( \IZ \oplus \frac{\IZ}{m \IZ} \right)^2, \qquad \Aut(B) \simeq B \rtimes \mu_n,\]
and the action of $\mu_n$ on $K$ is the one described in Remark \ref{rem: action}.
\end{thm}

\begin{rem}\label{rem: semidirect product identification}
The isomorphism $\Aut(X) / E \simeq K \rtimes \mu_n$ of Theorem \ref{thm: aut primary kodaira} sends the coset represented by the lift $F_{(v, \delta, \mu, \sigma)}$ to the element $\left( A \cdot \left( \begin{array}{c} \delta\\ \bar{v} \mu \end{array} \right), v \right)$, with $A \cdot \left( \begin{array}{c} \delta\\ \bar{v} \mu \end{array} \right) \in \Lambda_E^2 / (m\beta_2 \cdot \IZ)^2$. In this way, the multiplication is the standard one on a semidirect product, induced by the action described in Remark \ref{rem: action}.
\end{rem}

\section{Secondary Kodaira surfaces}\label{sect: secondary Kodaira}

Recall (see, e.g., \cite[Definizione 51]{Catanese}) that a \emph{secondary Kodaira surface} is a surface whose canonical bundle is non-trivial and which admits a finite unramified covering by a primary Kodaira surface. In this case, we can realize a secondary Kodaira surface as the quotient of a primary Kodaira surface by the action of a cyclic group $C_k$ of order $k$ (with $k = 2, 3, 4, 6$).

\begin{rem}
As the quotient of a primary Kodaira surface $X$ by the action of a finite subgroup of $E$ is again a primary Kodaira surface, we can assume that $C_k \cap E = \{ 1 \}$.
\end{rem}

We can be more precise: the canonical bundle of a secondary Kodaira surface $Y$ is torsion of order $k$, and its canonical cover is a primary Kodaira surface $X$. As a consequence, any automorphism of $Y$ lifts on $X$ and $\Aut(Y) = N_{\Aut(X)}(C_k) / C_k$.

Fix a lift
\[\gamma(z, \zeta) = (vz + \delta, \zeta + \mu z + \sigma)\]
for a generator of $C_k$. According to Theorem \ref{thm: aut primary kodaira} and Remark \ref{rem: semidirect product identification} its image in $K \rtimes \mu_n$ is the element
\[([\eta_\gamma], v), \qquad \text{with } \eta_\gamma = A \cdot \left( \begin{array}{c} \delta\\ \bar{v} \mu \end{array} \right) \in \Lambda_E^2.\]
As $C_k \cap E$ is trivial, the image $\tilde{C}_k$ of $C_k$ in $\Aut(X) / \Aut^0(X)$ is still cyclic of order $k$ and corresponds to the subgroup generated by $([\eta_\gamma], v)$.

We first compare the normalizer of $C_k$ in $\Aut(X)$ with the normalizer of its image $\tilde{C}_k$ in $\Aut(X) / \Aut^0(X)$.

\begin{prop}\label{prop: secondary kodaira short exact sequence 1}
Let $X \longrightarrow Y = X / C_k$ be the canonical cover of a secondary Kodaira surface $Y$. Then the normalizer of $C_k$ in $\Aut(X)$ fits into the short exact sequence
\[1 \longrightarrow E \longrightarrow N_{\Aut(X)}(C_k) \longrightarrow N_{\Aut(X) / \Aut^0(X)}(\tilde{C}_k) \longrightarrow 1.\]
\end{prop}
\begin{proof}
Since $E$ centralizes $C_k$ we have $E \leq N_{\Aut(X)}(C_k)$ and the short exact sequence
\[1 \longrightarrow E \longrightarrow N_{\Aut(X)}(C_k) \longrightarrow N_{\Aut(X)}(C_k) / E \longrightarrow 1,\]
and $N_{\Aut(X)}(C_k) / E$ is in a natural way a subgroup of $N_{\Aut(X) / \Aut^0(X)}(\tilde{C}_k)$. Focus on this last. As the order of $v$ is exactly $k$, we have that $N_{\Aut(X) / \Aut^0(X)}(\tilde{C}_k) = C_{\Aut(X) / \Aut^0(X)}(\tilde{C}_k)$. Let
\[f(z, \zeta) = (\varepsilon z + d, \zeta + p z + s)\]
be a lift for an element in $\Aut(X) / \Aut^0(X)$, then the element it represents centralizes $\tilde{C}_k$ if and only if $f \circ \gamma = e \circ \gamma \circ f$ for some translation $e(z, \zeta) = (z, \zeta + \tau)$. Equivalently
\begin{equation}\label{eq: centralizing condition}
\left\{ \begin{array}{l}
vd + \delta = \varepsilon \delta + d\\
\mu \varepsilon + p = vp + \mu\\
\tau + \mu d = p \delta
\end{array} \right. \longrightarrow \left\{ \begin{array}{l}
(1 - \varepsilon) \delta = (1 - v) d\\
(1 - \varepsilon) \mu = (1 - v) p\\
\tau = p\delta - \mu d.
\end{array} \right.
\end{equation}
Since $v \neq 1$, the first two equations determine $\left( \begin{array}{c} d\\p \end{array} \right)$ uniquely:
\begin{equation}\label{eq: d, p}
\left( \begin{array}{c}
d\\
p
\end{array} \right) = \frac{1 - \varepsilon}{1 - v} \left( \begin{array}{c}
\delta\\
\mu
\end{array} \right), \qquad \text{with} \qquad A \left( \begin{array}{c}
d\\
p
\end{array} \right) = \frac{1 - \varepsilon}{1 - v} A \left( \begin{array}{c}
\delta\\
\mu
\end{array} \right) \in \Lambda_E^2.
\end{equation}
The third equation in \eqref{eq: centralizing condition} then says that
\[\tau = \delta p - \mu d = \det \left( \begin{array}{cc}
\delta & d\\
\mu & p
\end{array} \right) = 0,\]
which means that an automorphism which normalizes $\gamma$ modulo $\Aut^0(X) = E$ already centralizes $\gamma$. In particular, it follows that
\[N_{\Aut(X)}(C_k) / E = N_{\Aut(X) / \Aut^0(X)}(\tilde{C}_k).\]
\end{proof}

Using the identification of Remark \ref{rem: semidirect product identification}, we can now describe the normalizer of $\tilde{C}_k$ entirely in terms of the semidirect product $K \rtimes \mu_n$.

\begin{prop}\label{prop: secondary kodaira short exact sequence 2}
Let $X \longrightarrow Y = X / C_k$ be the canonical cover of a secondary Kodaira surface $Y$. Let $([\eta_\gamma], v) \in K \rtimes \mu_n$ be the image of a generator of $C_k$ under the identification of Remark \ref{rem: semidirect product identification}. Then there is a short exact sequence
\[1 \longrightarrow \Gamma \longrightarrow N_{\Aut(X) / \Aut^0(X)}(\tilde{C}_k) \longrightarrow I_\gamma \longrightarrow 1,\]
where $\Gamma = \ker((1 - v): K \longrightarrow K)$ and $I_\gamma = \set{\varepsilon \in \mu_n \st (1 - \varepsilon) [\eta_\gamma] \in (1 - v) K}$.
\end{prop}
\begin{proof}
We regard $K$ as a $\IZ[\mu_n]$-module through the action described in Remark \ref{rem: action}. Let $(x, \varepsilon) \in K \rtimes \mu_n$. By construction, $([\eta_\gamma], v)$ is the image of a generator of $C_k$, hence both $([\eta_\gamma], v)$ and $v$ have the same order $k$. So we deduce that $(x, \varepsilon)$ normalizes $\langle ([\eta_\gamma], v) \rangle$ if and only if it centralizes it. Hence $(x, \varepsilon) \in N_{K \rtimes \mu_n}(\langle ([\eta_\gamma], v) \rangle)$ if and only if $(x, \varepsilon)([\eta_\gamma], v) = ([\eta_\gamma], v)(x, \varepsilon)$, i.e., if and only if $x + \varepsilon \cdot [\eta_\gamma] = [\eta_\gamma] + v \cdot x$. This is equivalent to
\[(1 - \varepsilon) \cdot [\eta_\gamma] = (1 - v) \cdot x,\]
hence the image of the projection $K \rtimes \mu_n \longrightarrow \mu_n$ is
\begin{equation}\label{eq: I_gamma}
I_\gamma = \left\{ \varepsilon \in \mu_n \st (1 - \varepsilon) [\eta_\gamma] \in (1 - v)K \right\}
\end{equation}
and its kernel is
\[\Gamma = \left\{ x \in K \st (1 - v) \cdot x = 0 \right\} = \ker((1 - v): K \longrightarrow K).\]
Thus we obtain the short exact sequence
\[1 \longrightarrow \Gamma \longrightarrow N_{\Aut(X) / \Aut^0(X)}(\tilde{C}_k) \longrightarrow I_\gamma \longrightarrow 1,\]
where $\Gamma$ is identified with the kernel of the above projection.
\end{proof}

We now pass from the canonical cover to the secondary Kodaira surface itself, and we begin with the following description of the group of connected components of $\Aut(Y)$.

\begin{prop}\label{prop: exact sequence secondary}
Let $X \longrightarrow Y = X / C_k$ be the canonical cover of a secondary Kodaira surface $Y$. Then $\Aut^0(Y) = E$ (the fibre of the Albanese map of $X$) and there is a short exact sequence
\begin{equation}\label{eq: ses secondary}
1 \longrightarrow \Gamma \longrightarrow \frac{\Aut(Y)}{\Aut^0(Y)} \longrightarrow \frac{I_\gamma}{\mu_k} \longrightarrow 1,
\end{equation}
\end{prop}
\begin{proof}
Taking the quotient of the exact sequence of Proposition \ref{prop: secondary kodaira short exact sequence 1} by the subsequence
\[1 \longrightarrow C_k \cap E = 1 \longrightarrow C_k \longrightarrow \tilde{C}_k \longrightarrow 1\]
we deduce the exact sequence
\[1 \longrightarrow E \longrightarrow \Aut(Y) \longrightarrow \frac{N_{\Aut(X) / \Aut^0(X)}(\tilde{C}_k)}{\tilde{C}_k} \longrightarrow 1,\]
and we deduce that $\Aut^0(Y) = E$ and $\frac{\Aut(Y)}{\Aut^0(Y)} = \frac{N_{\Aut(X) / \Aut^0(X)}(\tilde{C}_k)}{\tilde{C}_k}$.

Consider then $\tilde{C}_k \leq N_{\Aut(X) / \Aut^0(X)}(\tilde{C}_k)$: we have that
\[\Gamma \cap \tilde{C}_k = \{ 1 \}\]
and the image of $\tilde{C}_k$ in $I_\gamma$ is then $\mu_k \simeq \tilde{C}_k$. As a consequence, taking the quotient of the short exact sequence of Proposition \ref{prop: secondary kodaira short exact sequence 2} by the subsequence
\[1 \longrightarrow \Gamma \cap \tilde{C}_k = 1 \longrightarrow \tilde{C}_k \longrightarrow \mu_k \longrightarrow 1\]
we deduce the short exact sequence
\[1 \longrightarrow \Gamma \longrightarrow \frac{\Aut(Y)}{\Aut^0(Y)} \longrightarrow \frac{I_\gamma}{\mu_k} \longrightarrow 1.\]
\end{proof}

We conclude this section by providing an explicit description of the automorphism group of a secondary Kodaira surface. Observe that although the definition of $I_\gamma$ involves a choice of a lift of a generator of $C_k$, the resulting group is independent of these choices.

\begin{thm}\label{thm: list of groups}
Let $Y$ be a secondary Kodaira surface. Denote $X \longrightarrow Y = X / C_k$ its canonical cover and let $E$ be a fibre of the Albanese map of $X$. The automorphism group of $Y$ admits the following explicit description: $\Aut^0(Y) = E$ and there is a central extension
\[1 \longrightarrow E \longrightarrow \Aut(Y) \longrightarrow \Gamma \rtimes \frac{I_\gamma}{\mu_k} \longrightarrow 1.\]
The possible groups $\Gamma$ and $I_\gamma / \mu_k$ (which appear also in \eqref{eq: ses secondary}) are summarized in Table \ref{tab: tab1}, where $d_t = \gcd(t, m)$.
\begin{table}[ht]\centering
\begin{tabular}{|c|c||c||c|}
\hline
$n$ & $k$ & $\Gamma$ & $I_\gamma / \mu_k$\\
\hline
\hline
$2$ & $2$ & $(\IZ / d_2 \IZ)^2$ & $1$\\
\hline
$4$ & $2$ & $(\IZ / d_2 \IZ)^2$ & $1$ or $\IZ / 2\IZ$\\
$4$ & $4$ & $\IZ / d_2 \IZ$ & $1$\\
\hline
$6$ & $2$ & $(\IZ / d_2 \IZ)^2$ & $1$ or $\IZ / 3\IZ$\\
$6$ & $3$ & $\IZ / d_3 \IZ$ & $1$ or $\IZ / 2\IZ$\\
$6$ & $6$ & $1$ & $1$\\
\hline
\end{tabular}
\caption{The groups appearing in the exact sequence \ref{eq: ses secondary}.}
\label{tab: tab1}
\end{table}
\end{thm}
\begin{proof}
By Proposition \ref{prop: exact sequence secondary} it only remains to compute the groups $\Gamma$ and $I_\gamma / \mu_k$ explicitly.

Using the identification $K \simeq \frac{\Lambda_E^2}{((m \beta_2) \cdot \IZ)^2}$ provided by \eqref{eq: K} the group $\Gamma$ is identified with
\[\Gamma = \left\{ \left( \begin{array}{c}
y_1\\
y_2
\end{array} \right) \in \Lambda_E^2 \,\middle|\, (1 - V) \left( \begin{array}{c}
y_1\\
y_2
\end{array} \right) \in ((m\beta_2) \cdot \IZ)^2 \right\} \Big/ ((m\beta_2) \cdot \IZ)^2,\]
i.e.,
\[\Gamma = \left\{ \left( \begin{array}{c} a_1\\a_2 \end{array} \right) \in \IZ^2 \st (1 - V) \left( \begin{array}{c} a_1\\a_2 \end{array} \right) \in (m\IZ)^2 \right\} \Big/ (m\IZ)^2.\]
Using this description and \eqref{eq: I_gamma} we proceed with a case by case analysis on $n = 2, 4, 6$ and $k$ dividing $n$.

First, let $n = 2$ and $k = 2$. In this case we have $v = -1$ and so $V = -\id$, this means that $\left( \begin{array}{c} a_1\\ a_2 \end{array} \right) \in \Gamma$ if and only if $\left( \begin{array}{c} 2a_1\\ 2a_2 \end{array} \right) \in (m\IZ)^2$, hence we find that
\[\Gamma \simeq \left( \frac{\IZ}{d_2 \cdot \IZ} \right)^2, \qquad d_2 = \gcd(2, m).\]
As $n = k$, in this case we have that $\frac{I_\gamma}{\mu_k} = 1$.

Let $n = 4$ and $k = 2$. In this case we have $v = -1$, so
\[\Gamma \simeq \left( \frac{\IZ}{d_2 \cdot \IZ} \right)^2, \qquad d_2 = \gcd(2, m)\]
as before. The quotient $\frac{I_\gamma}{\mu_k}$ is a subgroup of $\frac{\mu_n}{\mu_k} = \frac{\IZ}{2\IZ}$, hence it can be either trivial or $\frac{\IZ}{2\IZ}$.

Let $n = 4$ and $k = 4$. In this case $v = \ii$ and we can choose a basis for $\Lambda_B$ such that $\Lambda_B \simeq \IZ \oplus \ii \IZ$. In this case then $V = \left( \begin{array}{cc} 0 & -1\\ 1 & 0 \end{array} \right)$ and so we have to solve
\[\left\{ \begin{array}{l}
a_1 + a_2 \equiv 0 \mod m\\
-a_1 + a_2 \equiv 0 \mod m
\end{array} \right. \longrightarrow \left\{ \begin{array}{l}
2 a_1 \equiv 0 \mod m\\
a_1 \equiv a_2 \mod m.
\end{array} \right.\]
As a consequence
\[\Gamma \simeq \frac{\IZ}{d_2 \cdot \IZ}, \qquad d_2 = \gcd(2, m)\]
and $\frac{I_\gamma}{\mu_k}$ is trivial as $n = k$.

Let $n = 6$ and $k = 2$, the same argument as before leads us to the conclusion that
\[\Gamma \simeq \left( \frac{\IZ}{d_2 \cdot \IZ} \right)^2, \qquad d_2 = \gcd(2, m).\]
The quotient $\frac{I_\gamma}{\mu_k}$ is a subgroup of $\frac{\mu_n}{\mu_k} = \frac{\IZ}{3\IZ}$, hence it can be either trivial or $\frac{\IZ}{3\IZ}$.

Let $n = 6$ and $k = 6$. We can choose a basis for $\Lambda_B$ such that $\Lambda_B \simeq \IZ \oplus \omega_6 \IZ$ with $\omega_6 = e^{\frac{1}{3}\pi\ii}$. Hence $v = \omega_6$, $ V = \left( \begin{array}{cc} 1 & -1\\1 & 0 \end{array} \right)$ and we have to solve
\[\left\{ \begin{array}{l}
a_2 \equiv 0 \mod m\\
-a_1 + a_2 \equiv 0 \mod m
\end{array} \right. \longrightarrow a_1 \equiv a_2 \equiv 0 \mod m.\]
So $\Gamma$ is trivial and, as $k = n$, also $\frac{I_\gamma}{\mu_k}$ is trivial.

Finally, let $n = 6$ and $k = 3$: as in the previous case, we have $v = \omega_6^2$, $V = \left( \begin{array}{cc} 0 & -1\\1 & -1 \end{array} \right)$. So we want to solve
\[\left\{ \begin{array}{l}
a_1 + a_2 \equiv 0 \mod m\\
-a_1 + 2a_2 \equiv 0 \mod m.
\end{array} \right. \longrightarrow \left\{ \begin{array}{l}
3a_2 \equiv 0 \mod m\\
a_1 \equiv -a_2 \mod m.
\end{array} \right.\]
As a consequence
\[\Gamma \simeq \frac{\IZ}{d_3 \cdot \IZ}, \qquad d_3 = \gcd(3, m)\]
and $\frac{I_\gamma}{\mu_k}$, being a subgroup of $\frac{\mu_n}{\mu_k} = \frac{\IZ}{2\IZ}$, can be either trivial or $\frac{\IZ}{2\IZ}$.
\end{proof}

\begin{rem}
Every group appearing in Table \ref{tab: tab1} is realized by some secondary Kodaira surface. For example, the maximal possibility $I_\gamma/\mu_k \simeq \mu_n/\mu_k$ is obtained by taking
\[\gamma(z, \zeta) = \left(vz, \zeta + \frac{\beta_1}{k} \right)\]
where $v \in \mu_n$ has order $k$. The cyclic group generated by $\gamma$ acts freely and properly discontinuously on $X$, so the quotient is a secondary Kodaira surface. In this case $\eta_\gamma = 0$, hence $I_\gamma = \mu_n$. The remaining possibilities are obtained by suitable choices of the generator $\gamma$, or equivalently of the class $[\eta_\gamma] \in K$. They can be realized as follows: we may assume that $\beta_3 = \beta_4 = 0$ and that $\alpha_3 = 1$ and $\alpha_4 = \ii$ in the case with $n = 4$ or $\alpha_4 = \omega_6 = e^{\frac{1}{3} \pi \ii}$ in the two cases with $n = 6$; then the following generators realize the desired groups
\begin{itemize}
\item for the case $(n, k) = (4, 2)$
\[\gamma(z, \zeta) = \left( -z - \frac{1}{2} \beta_1, \zeta - \frac{1}{2} \beta_1 z + \frac{1}{2} \beta_1 - \frac{1}{8} \beta_1^2 \right).\]
\item for the case $(n, k) = (6, 2)$
\[\gamma(z, \zeta) = \left( -z + \frac{\sqrt{3}}{3} \ii \omega_6 \beta_1, \zeta + \frac{\sqrt{3}}{3} \ii \omega_6^2 \beta_1 z + \frac{1}{2} \beta_1 - \frac{1}{6} \beta_1^2 \right).\]
\item for the case $(n, k) = (6, 3)$
\[\gamma(z, \zeta) = \left( \omega_6^2 z + \frac{\omega_6^2}{1 - \omega_6^2} \beta_1, \zeta + \frac{\omega_6^2}{1 - \omega_6^2} \beta_1 z + \frac{1}{3} \beta_1 + \frac{1}{9} \omega_6^2 (2 + \omega_6^2) \beta_1^2 \right).\]
\end{itemize}
In all three cases one has $\eta_\gamma = \left( \begin{array}{c} \beta_1 \\ 0 \end{array} \right)$, and \eqref{eq: I_gamma} gives $I_\gamma = \mu_k$.
\end{rem}

\section{Comparison with Fujimoto--Nakayama}\label{sect: comparison with FN}

Our computation begins with Kodaira's construction of primary Kodaira surfaces, but can be applied also with other constructions. For example, in \cite[$\S$6]{Fujimoto-Nakayama}, Fujimoto and Nakayama used this technique to find how a lift of a surjective endomorphism of a Kodaira surface looks like. Unfortunately, we found an inaccuracy in the proof of \cite[Proposition 6.4]{Fujimoto-Nakayama}. Here we give a correct statement of their result (we use the same notations and conventions of \cite{Fujimoto-Nakayama}); a proof can be given following the strategy exposed in the previous section.

Observe that the affine description obtained in Section \ref{sect: Kod surf} depends on the particular choice of the uniformization introduced by Kodaira. After passing to the coordinates used in \cite[$\S$6]{Fujimoto-Nakayama}, these affine transformations are expressed as polynomial automorphisms with a quadratic term in the second component. Thus Theorem \ref{thm: fn} should be viewed as the expression of Proposition \ref{prop: lift 2} in a different coordinate system.

\begin{thm}[{cf. \cite[Proposition 6.4]{Fujimoto-Nakayama}}]\label{thm: fn}
Let $f: X_{c, \delta} \longrightarrow X_{c, \delta}$ be a surjective endomorphism and let $h: T \longrightarrow T$ be the induced endomorphism with $\pi \circ f = h \circ \pi$. Suppose that
\[h_*: H^0(T, \Theta_T) \longrightarrow H^0(T, \Theta_T)\]
is the multiplication by $\alpha \in \IC$ with $\alpha L_\tau \subseteq L_\tau$. Then $f$ is induced from the automorphism
\begin{equation}
\begin{array}{rl}
\Phi(z, \zeta) = & \left( \alpha z + \beta, \abs{\alpha}^2 \zeta + \frac{1}{2} D_{\tau_B}(\alpha, 1) c \alpha z^2 + \right.\\
 & \left. + \pa{\sigma_{1, 0} + D_{\tau_B}(\alpha, 1) \pa{c\beta + \varepsilon - \frac{1}{2} D_{\tau_B}(\alpha, \tau_B) c}}z + v \right),
\end{array}
\end{equation}
of $\IC \times \IC$ where $\varepsilon = \delta - \frac{1}{2} c \tau_B$, $\alpha \in \IC \smallsetminus \set{0}$ is such that $|\alpha|^2 \in \IZ$ and $\beta, \sigma_{1, 0}, v \in \IC$.
\end{thm}

\begin{rem}
These are the consequences on other results in the literature, as far as the author knows. There is an issue in \cite[Corollary 3.3 and Corollary 3.4]{Shramov} (and in \cite[Theorem 1.1(iv)]{Shramov}), which was deduced as an immediate consequence of the particular form of the lifts provided by \cite[Proposition 6.4]{Fujimoto-Nakayama}. Anyway, the main result of that paper on primary Kodaira surfaces, i.e., \cite[Lemma 3.1]{Shramov} is correct. The same issue appears also in \cite[Lemma 2.4]{PS-BoundedAutomorphisms}.
\end{rem}

\bibliographystyle{alpha}
\bibliography{KodairaAuto}

\end{document}